# Excellence in Computer Simulation

by
Leo Kadanoff

The James Franck Institute
The University of Chicago
5640 S. Ellis Avenue
Chicago IL 60637 U.S.A.
E-mail: LeoP@Uchicago.edu

abstract

Excellent computer simulations are done for a purpose. The most valid purposes are to explore uncharted territory, to resolve a well-posed scientific or technical question, or to make a design choice. Stand-alone modeling can serve the first purpose. The other two goals need a full integration of the modeling effort into a scientific or engineering program.

Some excellent work, much of it related to the Department of Energy Laboratories, is reviewed. Some less happy stories are recounted.

In the past, some of the most impressive work has involved complexity and chaos. Prediction in a complex world requires a first principles understanding based upon the intersection of theory, experiment and simulation.

# The Best: Great Examples of Scientific Computing in the Heroic Period

I work within a Department of Energy supported research program at the University of Chicago, called the ASCI/Alliances Center for Astrophysical Thermonuclear Flashes. The main goal of ASCI, the alliance for scientific computing, is to gain knowledge and experience relevant for the construction of large-scale computer simulations, thus supporting computer research on complex systems, and thereby helping the DOE maintain our stockpile of nuclear weapons. My interest here is to provide an overview of the art and science of scientific computer simulation.

To begin: I am happy to tell you that the Scientific Laboratories of the Department of Energy have supported some of the very best scientific computing ever done. Indeed they may be said to have invented scientific computing. Below, I list some of the very best examples of scientific computing, and you will see a preponderance of Department of Energy supported work in this list.

In the "Monte Carlo" method a random number generator may be used to perform deterministic calculations[1]. The Rosenbluth's, the Teller's, Ulam, and Metropolis played major roles in putting the method together, and applying it in the "Metropolis algorithm" for calculating the properties of systems in thermodynamic equilibrium. When this calculation was performed, the idea that a numerical method could be built upon the fake, constructed randomness of the usual computer random number generator seemed almost magic. Indeed, in a famous incident, Geoffrey Householder said he would stop the computational work at Oak Ridge while he waited to be convinced that random number generators could work[2]. Today we have a fairly clear explanation of why the Monte Carlo method works, namely that for expectation values, the rapid increase in the number of configurations with energy is compensated by the rapid decrease of the Boltzmann factor, leaving a narrow window of states that actually have to be explored. Computers' random number generators are now frully accepted and reasonably well understood. From today's perspective, the most powerful and surprising feature that remains from this early calculation is the insight that the use of an inherently probabilistic algorithm can lead to an enormous compression in the number of computational steps required.

Monte Carlo calculations use a random number generator to simulate something else, often a system in thermodynamic equilibrium. In contrast, in a molecular dynamics approach the computer solves Newton's equations of motion and follows the trajectories of all the particles in the system. One of the first calculations of this kind was performed Fermi, Pasta, & Ulam[3] who studied a linear chain of atoms, coupled together with anharmomic forces, and driven by an initial large tweaking in one mode. A simulation was then used to follow the system as the 128 atoms jiggled about and transferred energy from one mode to another. After a time, the energy became satisfactorily mixed, as expected. However, they kept the computer going. And, at one magic moment, the system returned very nearly to its starting point, once again having almost all of its energy in one mode! These scientists had discovered, an almost integrable system, ....experimentally! By doing exploratory calculations in heretofore-unexplored territory, they stumbled across some wonderful new science.

Berni Alder and his collaborator Tom Wainright, working at Livermore National Laboratory, displayed an amazing mastery of the molecular dynamics method. They were involved in not one but two great discoveries. They looked at the motion of hard spheres bouncing off one another. To everyone's amazement, despite the purely repulsive interactions, they nonetheless saw a phase transition from a fluid state into a solid one.[4] Surprise number two is that these hard spheres, and indeed any colliding fluid particles, engender through their motion persisting correlations[5]. These "long time tails" remained a perplexing mystery for a long time, but now they are pretty well understood as a consequence of the hydrodynamic motion of the fluid as it flows past its own molecules.

Chaotic behavior, now characterized as sensitive dependence upon initial conditions, or "the butterfly effect", was discovered "accidentally" by Edward Lorenz[6] working with an early and very primitive program for solving linked sets of ordinary differential equations. It is said that one day the computer prematurely aborted a half-completed and recorded run. Lorenz went back and punched in the numbers corresponding (at least in their first few digits) to the results used by the computer on the previous day. To his amazement and our subsequent edification, the resulting run was entirely different from the one of the previous day. The high order digits were seen to matter, a lot. Chaos was discovered!

My last "great" example comes from Los Alamos[7]. Mitchell Feigenbaum was using a not-very-fancy desk calculator to study the properties of models, previously investigated by Ulam, based upon a quadratic formula. The calculator took a number, squared it, formed a linear combination from the square and some fixed coefficients, and generated a new number. That simple process was carried on through many steps by Feigenbaum and his trusty computer. And patterns emerged! Wonderful, unexpected patterns-- showing how systems even simpler than Lorenz's could become "just a little bit" chaotic. And a exciting little world opened up, unexpectedly.

In all these examples, the scientists involved discovered and explored entirely new pieces of science. In the Monte Carlo case, the novelty was in the use of a new kind of computer algorithm which new conceptualizations of possible computations.   In the other cases, the use of a highly simplified model combined with a new kind of hardware or a new kind of calculational technique permitted the first scientific investigation of new domains of physics.  These researchers discovered new scales of length or time: long term recurrence or long-range order.   They found and were receptive to new scientific ideas, which were general and applicable to a broad range of systems.  Subsequent experiment, theory, and simulation has studies each of these ideas in great detail, and taken them much further.

## Good  Recent  Examples.

We may suspect that the heroic age is now passed. The nature of discoveries is now somewhat different.  I reach for recent examples to describe the best things that computational people are now doing.

One of the best pieces of science done in recent years is the discovery of neutrino mass and neutrino oscillations. The first hint of this major discovery came from a discrepancy between the flux of neutrinos from the sun measured by Ray Davis and others compared with the flux predicted by the computer models of solar reactions and activity.  In order for the discrepancy to be taken seriously, one had to believe in the accuracy and reliability of the solar models[8]. It was persuasive that the experiments were done by extremely competent people who believe in their results.  Another persuasive factor was an observational program seeking additional tests for the models.  The models were verified by a comparison with seismographic data recording wave activity within the sun.  The seismographic predictions of the models fit the observations.  Now with increased credibility of the models the original discrepancy in neutrino fluxes was seen to be a serious problem. Something had to give.   Eventually the part of the picture concerned with elementary particle physics had to be modified. The then-accepted theory assumed that neutrinos have zero mass. This idea was abandoned to fit the solar data.  Later observations have supported this change.

The big discovery was made by the experimentalists who observed and counted  neutrinos.  Major credit also went to the theorists who held the scientific program together, particularly John Bahcall.  The computer model-builders were a third, and quite essential,  part of the enterprise.   All together, these investigations have produced a major unexpected advance in our understanding of the fundamentals of the universe.

Another interesting recent example involves the problem called single bubble sonoluminescence. This phrase describes a situation in which a resonant acoustic field forms a bubble and excites it strongly enough so that the bubble emits light and becomes visible[9].   Here, once again, the major discoveries were made by the experimentalists. It was they who discovered the phenomenon and started a major program of work in the area. By 1997 they had set out the values of some of the most important parameters describing the system. As this understanding of the subject developed, simulations were performed and closely tied to the theoretical and experimental work.  Thus, very many of the proposed steps and advances were looked at, checked, and more deeply understood because of the simulations. The physical and chemical processes involved are quite complex, and the simulations enabled the workers to tie everything up into one bundle. In this example simulation played a major integrating role, but never led the advances in the field.  I would judge that this is the common role for simulation in "table-top" science.  I expect that this example will serve as a paradigm for future investigations in that a scientific and technical problem involving many branches of science has been largely understood through the interlinked efforts of many investigators who often used simulations to check their arguments.

My third example is drawn from computer simulations of the cosmology of the early universe[10]. Our colleagues working in this area have a hard task. Experiments are impossible.  They can only observe-- but never manipulate-- their system.  They are much further away from their system, in both space and time, than are the solar physicists, so even their observational data are quite limited in scope.  For this reason, they must be more dependent upon simulations than people in most other branches of science.  They construct entire universes, intended to be reasonably realistic, within their computers. They start from dark matter and baryons and then observe the model bring together clusters on a variety of scales.  Step by step the computers make objects on many scales,  down to the size of galaxies.  Can their constructed universes give them real insight into the processes involved?  One might worry that the model-making gives too much freedom, so that simulators will always be able to fit the known facts with a wide variety of schemes.   The simulators disagree. Bertschinger states that the "main use [of the calculations] has been and continues to be the testing of the viability of cosmological models of structure formation." [ibid page 599] The work takes the theoretical conceptions of the field, cast them into the form of specific models and then runs them. Many models simply blow up, yielding nothing sensible.  The remaining ones give well-defined results  which can be analyzed to see whether they agree with the observations. Many models fail at this stage.

In the long run the simulators hope to strain out all but the one, correct physical model for the development process.  This filtering is a brave goal, which the participants believe they can achieve.  I cannot tell whether they are being too optimistic.   At the moment, several different models apparently work quite well: "Recent high-resolutions simulations compare remarkably well with many aspects of the observed galaxy distribution." [ibid page 632]

In all three of these recent examples, the role of simulation was to work with theory, observation, and experiment to serve as a cross-check on the other modes and thereby to increase the confidence of the investigators that they understood phenomena that were not otherwise open to observation.  In the solar neutrino example, the simulations made the investigators confident that they understood the solar behavior and thus they were able to able to locate an error in our assumptions about neutrinos. In the sonoluminescence case, the simulations were a necessary part of putting together an intricate puzzle.  The early universe investigators hope and expect to weed out incorrect mechanisms and theories by carefully testing their consequences and comparing with our observations about the universe. In each case, the simulation works by being part of a carefully constructed program of activity.

## Not so Good: Optimization of Enthusiasm/Misjudgment

In a recent development, a provocative and controversial experiment conducted at Oak Ridge suggested that fusion was occurring in deuterated acetone--via a process involving resonance excitation of bubbles.   The reporting paper[11] involved both experimental work and computer simulations. "[A] roughly ten-fold increase in the external driving pressure was used in the calculations" beyond the pressure directly produced by the experimental situation "to approximately account for the effect of pressure intensification within the imploding bubble clusters".   As a result their "[h]ydrodynamic shock code simulation supported the observed data". It is remarkable that the refereeing process for a high visibility paper allowed an apparently uncontrolled approximation in a key step in the computer calculation.   Subsequent work[12] seemed to disprove the Oak Ridge experimental result. But that is not the point here.  Because of the "roughly ten-fold increase", the simulation was sufficiently uncontrolled so that it neither supported nor could refute the experiment.  It was simply beside the point. Neither the authors nor the editors should have permitted it to be published.

This example makes one ask what kind of quality control is appropriate for a computer calculation used to check a provocative experimental result.

This problem is broader than this one paper. The whole early history of single bubble sonoluminescence[13] required step by step work to eliminate provocative but incorrect mechanisms. A set of early experiments by Barber et. al[14] reported very short widths for the emitted pulse of light. This short width then opened the door to novel mechanisms for explaining the total intensity of the emitted light. Later developments suggested that the short pulse width was a misstep by the experimentalists. In contrast to the excellent work in sonoluminescence in the post-1997 period, reported above, this misstep led the simulators and theorists quite astray. A host of incorrect speculations and mechanisms ran through the field, intended to explain the "observed" behavior. Despite one essentially correct simulation[15], the pre-1997 simulations[16] did almost nothing to weed out these incorrect discussions, undercutting ones hope that simulations might provide a good tool for such weeding. (Note that this weeding was one of the main goals of the astrophysical modeling.) Instead the speculations continued unhindered until an experiment by Gomph and coworkers[17] showed that the pulse width was much longer than previously believed. This implied a lower temperature for the emitting drop. After this, attention turned away from the incorrect mechanisms so that-- as reported above-- theory, experiment, and simulation began to produce a consensus about what was going on.

The examples of the Oak Ridge paper and some of the earlier sonoluminescence simulations suggest that the models might have been directed toward the wrong goals. Apparently, rather than being used for a process of checking, criticism, and elimination of incorrect possibilities they were often used to support and exemplify the presumptions of the scientists involved. A program of modeling should either elucidate new processes or identify wrong directions. Otherwise, there is no point in carrying it out.

Another example, which might be entirely mythical, involves a transportation investment model said to have been put together in Britain with the goal of getting the best transportation system while minimizing public spending. The model involved a broad mix of roads, rail, and public improvements. The goal was overall maximization of benefits, taking into account public spending and also the value of time saved. All costs and benefits were converted into pounds and an overall optimization was sought and achieved.

The next step was to bring in an outside group of experts to study the model's recommendations and to bring together a plan for implementing them. This group noticed several apparent anomalies. The strangest, according to the story, was the elimination of all spending for improving pedestrian street-crossings. This result was considered peculiar, especially since the value of pedestrian time saved was included in the model. A careful look explained how the conclusion was reached. The decreased spending had the effect of increasing accidents at the crossings. According to experience, and also the model, the major result would be increased deaths among older pedestrians. Thus spending on pensions would be reduced. The model counted this outcome as a benefit.

The government that had paid for the modeling was not amused[18].

This outcome brings us to a moral: The transportation study failed because the modeling had been done too mechanically, without enough thinking about either the actual processes going on or the actual goals of the sponsors. They did not realize that the design goals were actually multidimensional. They did not ask why did we get this outcome? Modeling efforts should include theory and common sense. The two examples relating to bubbles have a different moral. In these cases the simulations, both the several wrong ones and even the essentially correct simulation of Vuong and Szeri (reference 15), did not effectively refute the incorrect experiments. Instead the simulations were effectively trumped by experiments, which the community judged to be decisive.

## Present Challenges: I. Convective Turbulence

In this section and the next I describe some work involving simulations in which I have played some role. In Rayleigh Bénard flow, one observes a fluid in a box heated from below and cooled from above. One describes this flow in terms of a parameter, called the Rayleigh number that gives a dimensionless measure of the strength of the heating. The higher the Rayleigh number, the more turbulent is the system. To compare with other turbulent systems, I might say that the Rayleigh number is roughly the square of the Reynolds number or the fourth power of the Taylor Reynolds number.

A little heating of system from below causes no motion of fluid. With increased heating and increased Rayleigh numbers one sees, first motion, and then chaos. At Rayleigh numbers above roughly $10^8$, turbulent flows and structures are formed, as seen for example in Figure 1. As shown in the cartoon of Figure 2., the heated box contains many structures including plumes, waves, and jets. How far are we from examining this experimental behavior in computer simulations?

There exist good simulations, in both two and three dimensions, but the three-dimensional simulations don't resolve the structures seen in the experiments. Experiments now reach to Rayleigh numbers as high as $10^{19}$. Simulations hardly go beyond $10^{12}$, because of limitations caused by resolution and computer time.  Theory suggests phase transitions, qualitative changes in behavior, at roughly $10^8$, $10^{11}$, and $10^{19}$.  Theorists are unsure of what will happen and in fact consider a very large range of possibilities.   Simulations cannot hope to reach directly into the domains touched by theory and experiment. Nonetheless, we are beginning to learn how to use theoretical ideas to extrapolate simulation results from lower Rayleigh numbers to higher ones.  The simulations provide detailed information to help us see what is really happening, in much more detail than the experiments can now provide.  The high Rayleigh number experiments' data-generation is limited by the design, manufacture,  and placement of delicate and tiny temperature measuring devices.

One recent example is a simulation done by Marcus Brüggen and Christian Kalser (Nature **418** 301 (2002) describing a hot bubble toward the center of a galaxy as in see Figure 3.  Because we cannot see into the galactic center, this bubble can only be "observed" through computer simulation.  Nonetheless the authors are confident that they have caught some of the essential features of heat transfer in this region.

Table 1 provides a comparison of what we might gain from experiment and what we might gain from simulation.  Clearly both are necessary.  Theory is also required to extrapolate the simulation-result into a physically interesting situations. More broadly, theory is required for simulators

> to assess reliability of algorithms
> to make better algorithms
> to help define what's worth "measuring"

Theorists also often help bring it all together-- recall the work of Oppenheimer, Teller, and Bahcall.  Of course in the ideal world you would have a scientist who could do it all, like  Leonardo do Vinci or Enrico Fermi.  But usually, different people have different quite specialized skills. To solve hard problems one must make all the kinds of scientific skills work together and, in the end, pull in the same direction.

| quantity | simulation | experiment |
|---|---|---|
| turnovers | five or ten | thousands |
| Ra | up to $10^{11}$ | Up to $10^{14}$ |
| runs | Few and costly | many |
| flexibility | low | high |
| measure | anything | Very few things |
| precision | Often very high | variable |
| equations | well known | Often unknown |
| Small variation in initial data | easy | impossible |

Table 1. Experiment and Simulation Complement one Another. In the first four rows experiment can do better because it runs longer, with more extreme flows, more repetitions, and hence more flexibility. But experimentalists can measure few things, relatively imprecisely, in hard-to-control situations. They also cannot change the initial data just a little and run again. So both approaches are necessary.

# Present Challenges: II. Jets and Sprays

Here we observe dielectric and conducting fluids moved by an electric field.  Let us start from an experiment.

Experimental work done by Lene Oddershede and Sidney Nagel start from oil floating on water as depicted in Figure 4. They apply a strong electric field, with the regions of strongest field strength being near the curved electrode sitting in the oil. The lower fluid, the one with the higher dielectric constant, is pulled upward toward the stronger electric field.   Thus, in the first few panels, we see that the water forms itself into a bump.

Here is a nice simple problem which we might perhaps use as an exercise in a simulational partial differential equations course.  The flow looks simple and easy to understand.  But, in the real world, surprises are possible, even likely (see the last few panels).  After a time, the water bump forms itself into a sharp point.  Then, starting from the point, something starts moving through the oil.  In the next to last frame, that motion resolves itself int a jet of charged fluid.  Then in our final frame, the fluid breaks up into many tiny droplets.

Complex systems sometimes show qualitative changes in their behavior. (Here a bump has turned into lightning and rain.) Our simple problem has developed new phenomena and new scales

Experiment is very good at finding unexpected behavior and describing its overall characteristics. Theory can often explain what is going on.  Then, after an appropriate pause for algorithm development,  simulations can test the ideas and fill in the details.

More recently, my student Moses Hohman has established the basic mechanism for the production of rain by doing a simulation investigating the linear stability (or rather instability) analysis of a charged jet.    Working with Michael Brenner, M. Shin, and G. C.Rutledge he looked for and saw a whipping instability in the motion.  This instability produces a turning motion, rather like that of a corkscrew.  The drops are presumed to be thrown off by the spinning.

In a parallel effort, my student Cheng Yang has looked at the process of singularity formation in the interface between two unlike dielectric fluids in motion in an electric field.  He was looking for the structure formed very near the singularity. He found a surprise, a result contrary to our initial presupposition.  From the previous literature, especially the work of G.I. Taylor,  I expected to see the formation of a static cone-like structure which could have a persistent existence in the electric field. Yang actually found (Figure 5) a transient dynamical conical structure, which formed for an instant and then broke up. As his thesis adviser, I am more than slightly proud that his simulation found something unexpected, and that he stuck to his guns long enough to convince his thesis committee that his result was both surprising and correct. So often, simulations only yield what was desired from the beginning.

## Present Challenges: III. The Rayleigh Taylor Instability

This instability has been an important focus of recent work, especially within the ASCI program. The instability can arise whenever a heavier fluid sits on top of a lighter one. If the interface between the two remains horizontal, nothing happens. However, a wrinkling of the surface can produce a cascade of changes in which jets and plumes of the heavier fluid penetrate into the lighter one and vice versa.

Some experimental studies of this situation have been carried out. In Early in the ASCI programs, for administrative reasons, a decision was made to have the program concentrate upon simulations-- with only minor input from experiment. More recently, apparently, the weakness in this unbalanced approach was recognized, resulting in an increased emphasis upon experiment. Some of the Rayleigh Taylor work, however, was affected by the earlier, unbalanced style.

Many important simulations of the Rayleigh Taylor system have been performed. In order to see the fully-developed instability some major simplifications of the physical model are required. Since ASCI is interested in large Reynolds numbers, the viscosity is usually neglected in the simulations. Further, to maximize the effect of the instability one neglects the surface tension in the interface between the fluids. These choices have been made to speed up the simulation. They do that. However, the problem which remains is technically "ill-posed" in that one cannot prove that it is mathematically meaningful. The practical meaning is that one cannot promise that different approximation approaches will converge to the same answer.

The outcome has been, to say the least, quite interesting. A group of studies have been put together-- all aiming to measure the degree of penetration of one fluid into another[19]. The penetration is determined in terms of a coefficient called $\alpha$ which measures the extent of the mixing zone relative to a purely ballistic motion of the fluids. An experiment measuring this quantity has been compared to half a dozen different group's simulations, all starting from identical initial conditions. The results fall into two groups. The experiment[20], the theory[21], and one of the simulations[22] show an $\alpha$-value of roughly 0.06; the other simulations give $\alpha$ of the order of 0.03 or less. (See Figure 6) Another study (reference 19) takes a different tack by looking at a single penetrating region by using periodic boundary conditions. (See Figure 7.) Note that the flow is extremely complex and quite sensitively dependent upon the computational resolution. If one takes the generated pictures at their face value one would conclude that the shape of the interface given by the simulation will never converge. On the other hand, there is some indication of convergence of the value of $\alpha$. We still do not know if the approximation of zero surface tension and viscosity make any sense, and if the value of $\alpha$ obtained in this way is meaningful.

To drive this point home, we look at one more example. Figure 8 shows four calculations of the mixing of the spray produced by a breaking wave. All four describe the same, ill-posed problem: wave motion without surface tension or viscosity. All four start from the same initial data. All four have the same value of the "wind" driving the wave. The only differences are in calculational resolution. And in the answers. The patterns of spray look quite different. The graph on the far right shows that the amount of mixing is not only quite different, but is a non-monotonic function of resolution. In short, much more work will be required before one can, with full reliability, estimate the mixing from this calculational method.

The problems with these calculations point out, once more, the well known fact that finding a valid answer from a computer simulation can be a matter of some subtlety. For example, the calculation shown in Figure 3 has a range of validity which must be regarded as unknown because the numerical method must still be regarded as unproven. The calculation describes events at the center of a galaxy. We are interested in having an accurate picture of what goes on here, but we can well afford to wait for the further scientific developments which can be expected to tell us more about the accuracy of the calculational method. In other cases, however, we may need accurate answers to questions involving highly turbulent flows. Unfortunately, we have no proven way of getting them.

## Conclusions

I. To maintain a national capacity for understanding the development of complexity and multi-scale phenomena, we should support first principles studies of a variety of different complex systems. Each such study requires a balanced and interdisciplinary program of research in which theory, simulation and experiment work together to ask and answer incisively posed questions.

II. The goal of my group's research at Chicago is to ask important questions about the world. We solve simple model problems, like the two that I have discussed here, and then ask questions like:

> How does complexity arise?  Why is chaos often observed?
> What dramatic events occur in the fluid?   Are they commonplace?
> Why do fluids naturally form structures?

A parallel goal is to teach students to ask incisive questions.

These are good problems for students because they are small enough so that they can be solved quickly.  They are also down-to-earth enough so that each and every student can appreciate what they are about.

III. In the world outside the schools, we simulators have an important role to play in as a part of the teams within scientific and engineering groups devoted to understanding and to design and development.  In the past, we have sometimes appeared in a supporting role, filling in the details in understandings constructed by others.  We may wish to be more incisive, pointing out where the design won't work, the theory won't hold water, the experiment is wrongly interpreted.  We may wish to be more creative, using our simulations to point the way to the overall understanding or the good design. Then we can expect our work to be evaluated and tested by the hands-on and the pencil-and-paper people who will also form a part of our scientific and engineering world.  Such an give-and-take approach forms the basis of good design and understanding.

IV. Conversely, if our work only justifies and explains the work done by designers and experimentalists, if we simulators never say that the other guys are dead wrong, then we deserve a situation in which simulation is relegated to the position of a third and lesser branch of science, way behind either experiment or theory.

## Acknowledgements

This work was supported in part by the Department of Energy through the ASCI/FLASH program and by the University of Chicago MRSEC and via a NSF-DMR grant.  I would like to thank Alexandros Alexakis, Michael Brenner, Alan Calder, Sascha Hilgenfeldt, Robert Laughlin, Steve Libby, Detlef Lohse, Sidney Nagel, Robert Rosner,  Andrew J. Szeri, Penger Tong, Kequing Xia, and Yuan-Nan Young for helpful comments and for help in getting together the figures.

# Figures

Figure 1. A shadowgraph showing the spatial distribution of thermal plumes in a Rayleigh Bénard cell.  The fluid is dipropylene glycol, which has a rather high viscosity (Prandtl number =596) so that the pattern of plumes appears in a simplified form.   The Rayleigh number is $6.8x10^8$. The picture was taken by S. Lam in Kequing Xia's laboratory and will appear in a joint publication with Penger Tong.

Figure 2.  Cartoon of Box. This shows the central region of the box, a mixing zone containing  plumes, and very thin boundary layers at top and bottom.  The plumes are believed to arise from the spray thrown up by waves travelling across the boundary layer.  This cartoon is redrawn from X.-L. Qiu and P. Tong, Phys. Rev. E, vol. 66, 026308 (2002).

Figure 3 A simulation done by Marcus Brüggen and Christian Kalser (Nature **418** 301 (2002) describing a hot bubble rising from the center of a galaxy. In these pictures, gravity points leftward.  The upper picture gives an earlier time; the lower one shows a later snapshot of the thermal bubble. The simulations describe two-dimensional flow with zero viscosity and zero thermal conductivity The color coding describes density.

Figure 4. Shows the experimental production of a singularity at the onset of an electrohydrodynamic spout ( Lene Oddershede and Sidney R. Nagel, Phys. Rev. Lett. vol. 85, 1234-1237 (2000)).  Two fluids, "oil" above and "water" below,  are separated by an interface. There is a strong electric field pointing into the electrode shown at the top of each frame.  This field carries the water, with its higher dielectric constant, upward into the region of strong electric field. Eventually, the interface comes to a point, breaks down, a discharge is produced, and generates many small droplets of water in the oil.

Figure 5.  Shows a computer simulation of the production of a singularity in a situation in which two fluids with different dielectric constants are separated by an interface.  The electric field generates polarization which then produces forces on the surface of the drop.  Surface tension provides additional forces.  The first frame shows the initial and final shapes of the drop.  Originally we have an ellipsoidal shape; after a time the drop develops cone-like points on the ends.  The second frame shows how the cone gradually sharpens.  The last shows that there is indeed a singularity in that the velocity of the tip diverges at a critical time.  This picture was taken from the Ph.D. thesis of Cheng Yang at the University of Chicago.

Figure 6   The Rayleigh Taylor instability.  The initial state was a gently wavy interface separating a high density fluid from a low density one. Gravity (pointing up!)  then destabilizes the interface, producing the mixed regions shown. Unmixed regions are transparent. Red yellow and green show successively higher densities. This simulation assumes that both viscosity and surface tension are negligibly small.  Taken from reference 19.

Figure 7 The Rayleigh Taylor instability once more. This is an unpublished picture given to me by Alan Calder.  A somewhat similar picture appears in reference 19.  This calculation is done with a simple "one-bump" initial state. The effect of resolution is studied by employing resolutions differing by a factor of two in successive panels.  Note that the results change considerably with resolution. The highest resolution picture is qualitatively different from the others in that the left-right symmetry is broken.

Figure 8 Wave breaking at a White Dwarf Surface.  This figure shows the result of a wind-driven instability on the surface of a star. Surface tension and viscosity are assumed to be negligibly small.  The different panels show, once more, resolutions differing by a factor of two with the same initial condition and at the same time.  On the right, one sees plots of mixing versus time for these different resolutions.  The take-home message is that resolution matters both in the profile and also in the mixing.    Alexandros Alexakis did the simulations for these images. The images were supplied to me by A. Calder.   Similar results appear in A. C. Calder, A. Alexakis, L. J. Dursi, R. Rosner, J. W. Truran, B. Fryxell P. Ricker, M. Zingale, K. Olson, F. X. Timmes, P. MacNeice "Mixing by Non-linear Gravity Wave Breaking on a White Dwarf Surface" Proceedings of the International Conference on Classical Nova Explosions, Sitges, Spain, 20-24 May 2002 (AIP Conf.Proc. 637 (2003) 134-138.)

---

[1]. Metropolis, N., A.W. Rosenbluth, M.N. Rosenbluth, A.H. Teller, and E. Teller, 1953; Equation of State Calculations by Fast Computing Machines, Jour. Chemical Physics, V. 21, No. 6, pp. 1087 - 1092. This paper was cited in Computing in Science and Engineering as being among the top 10 algorithms having the "greatest influence on the development and practice of science and engineering in the 20th century."

neglected. These codes underestimated the damping mechanisms and hence produced a very strong shock, which would, in the approximations used by the investigators, produce an infinitely high temperature. (See also the parallel work by W. C. Moss, et. al. ("Hydrodynamic simulations of bubble collapse and picosecond sonoluminescence" Physics of Fluids **6,** 2979-2985, (1994)). Later simulations by Vuong and Szeri (see previous reference ) cast doubt upon the relevance of shocks to the observed behavior of sonoluminescence. However, the field did not turn around until new experimental results caught people's attention.

# Figure 1

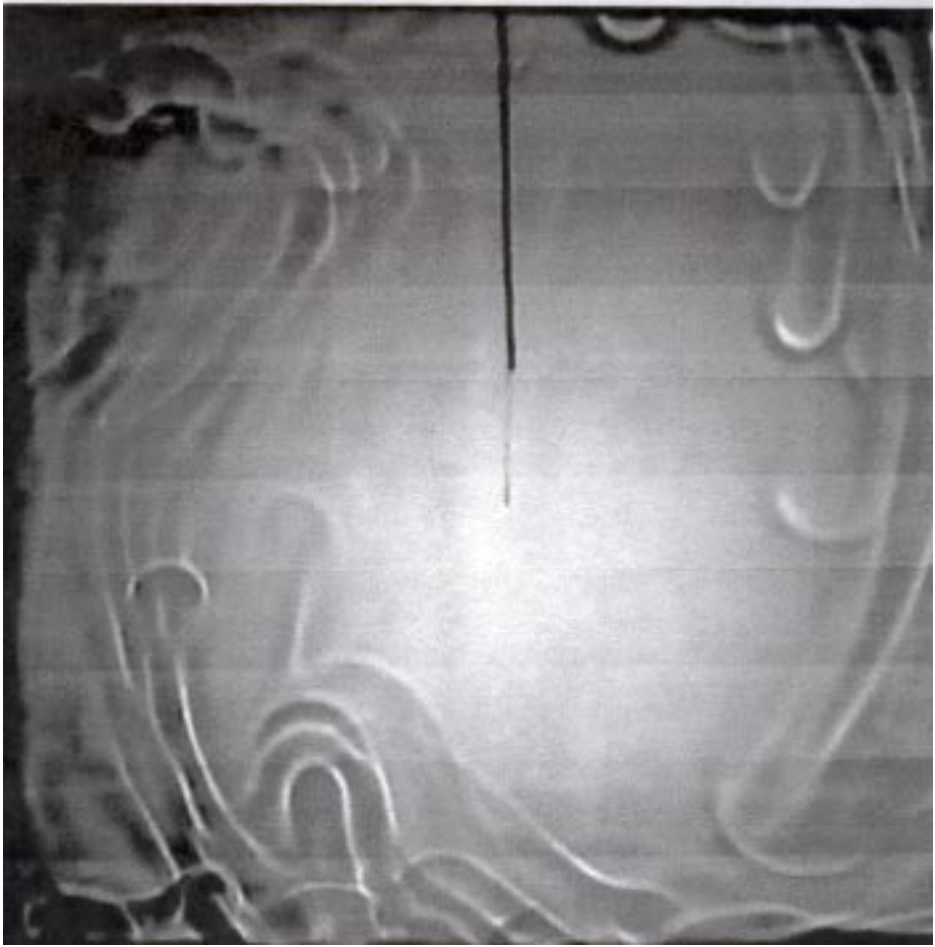

# Figure 2

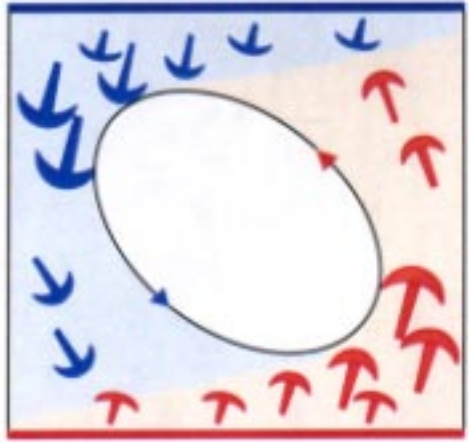

# Figure 3

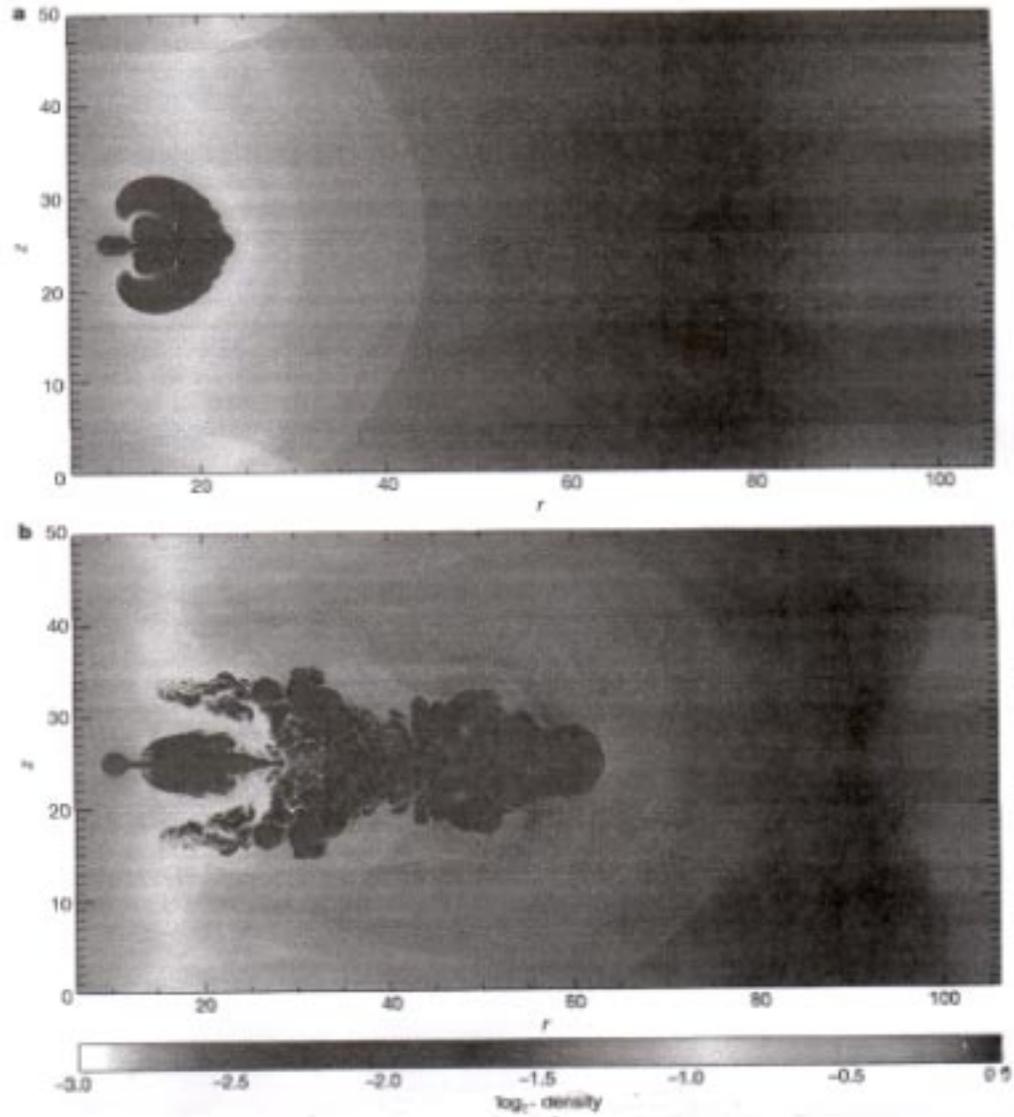

# Figure 4

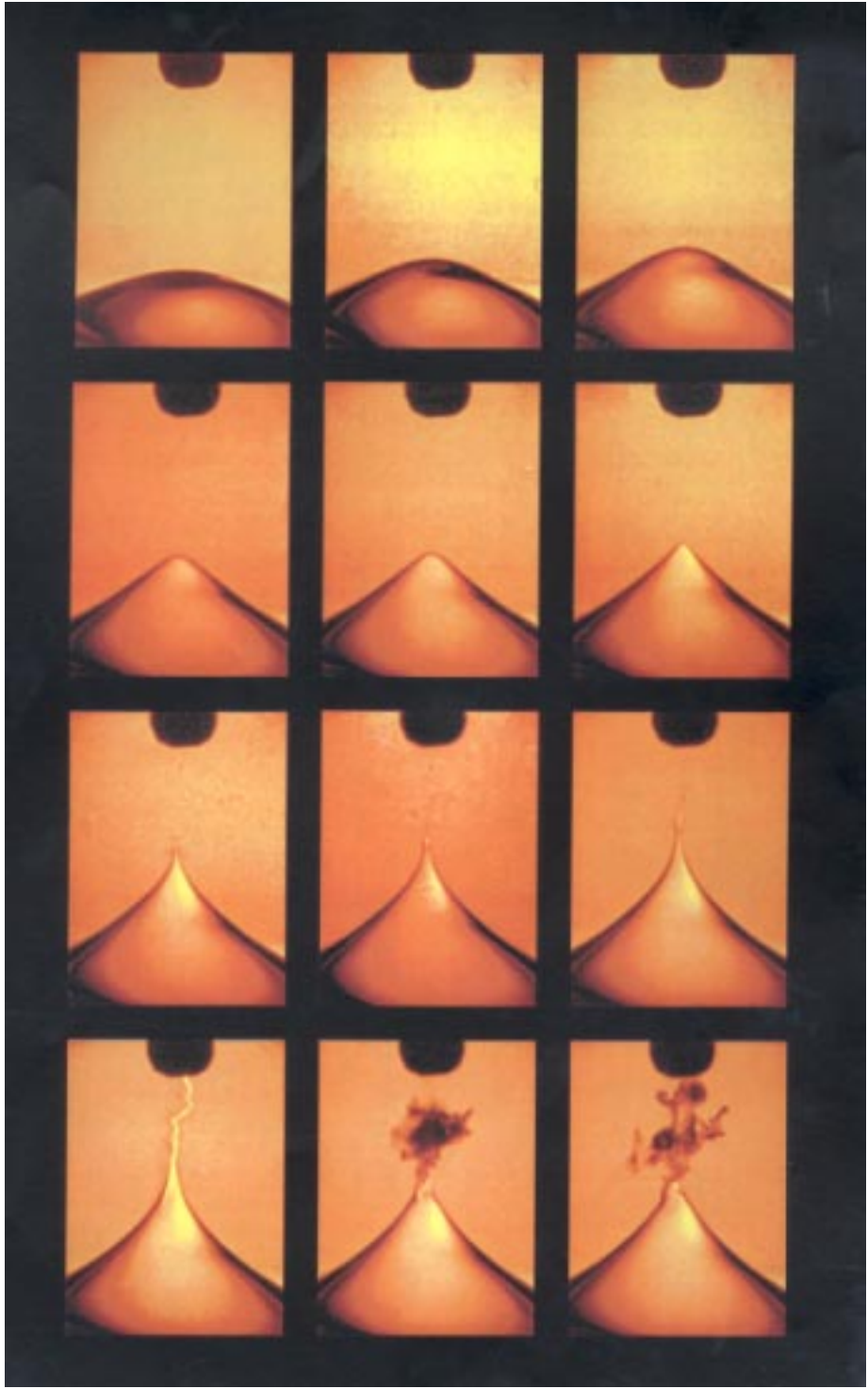

Figure ⟶

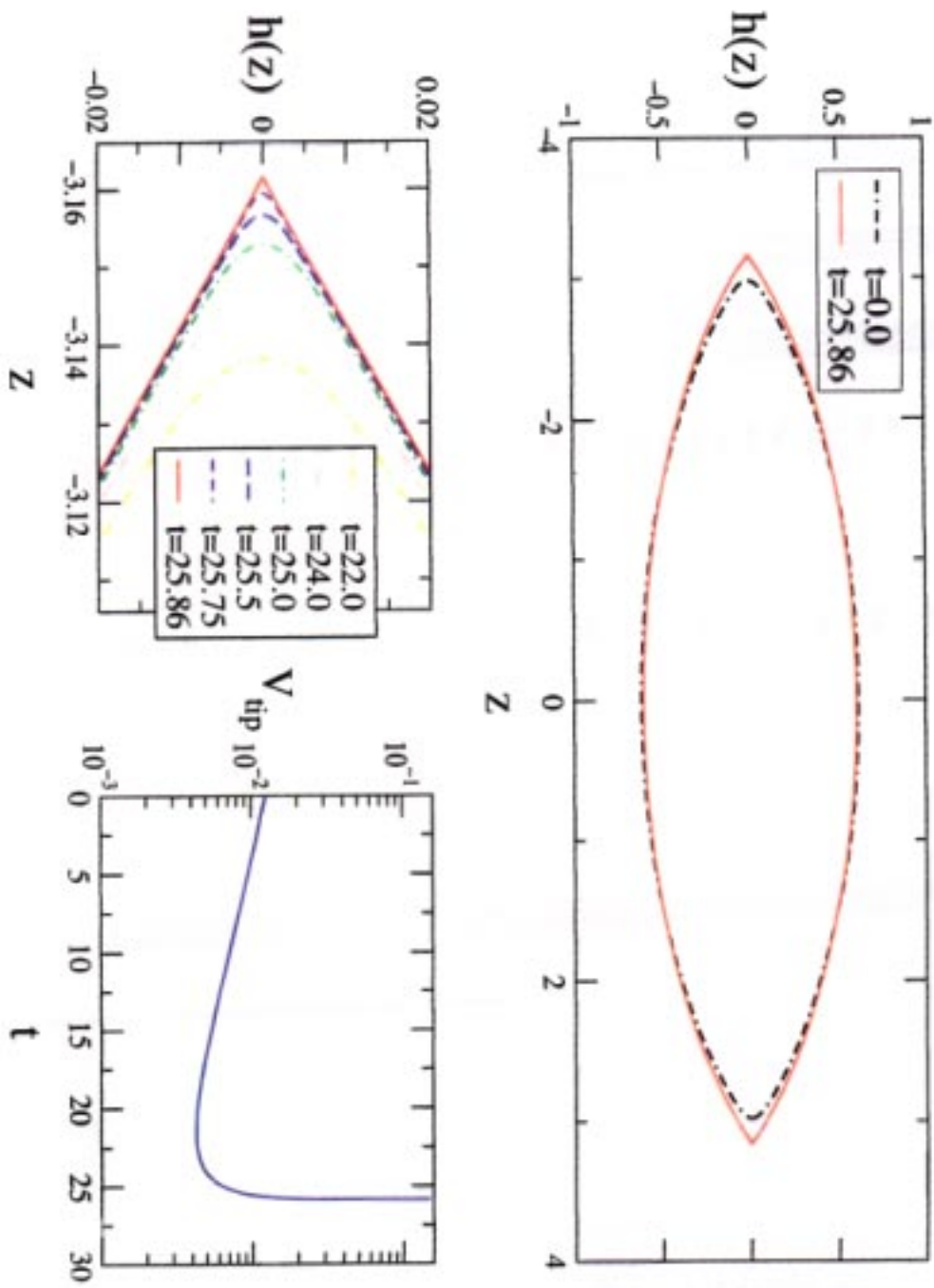

Figure 6 ↑

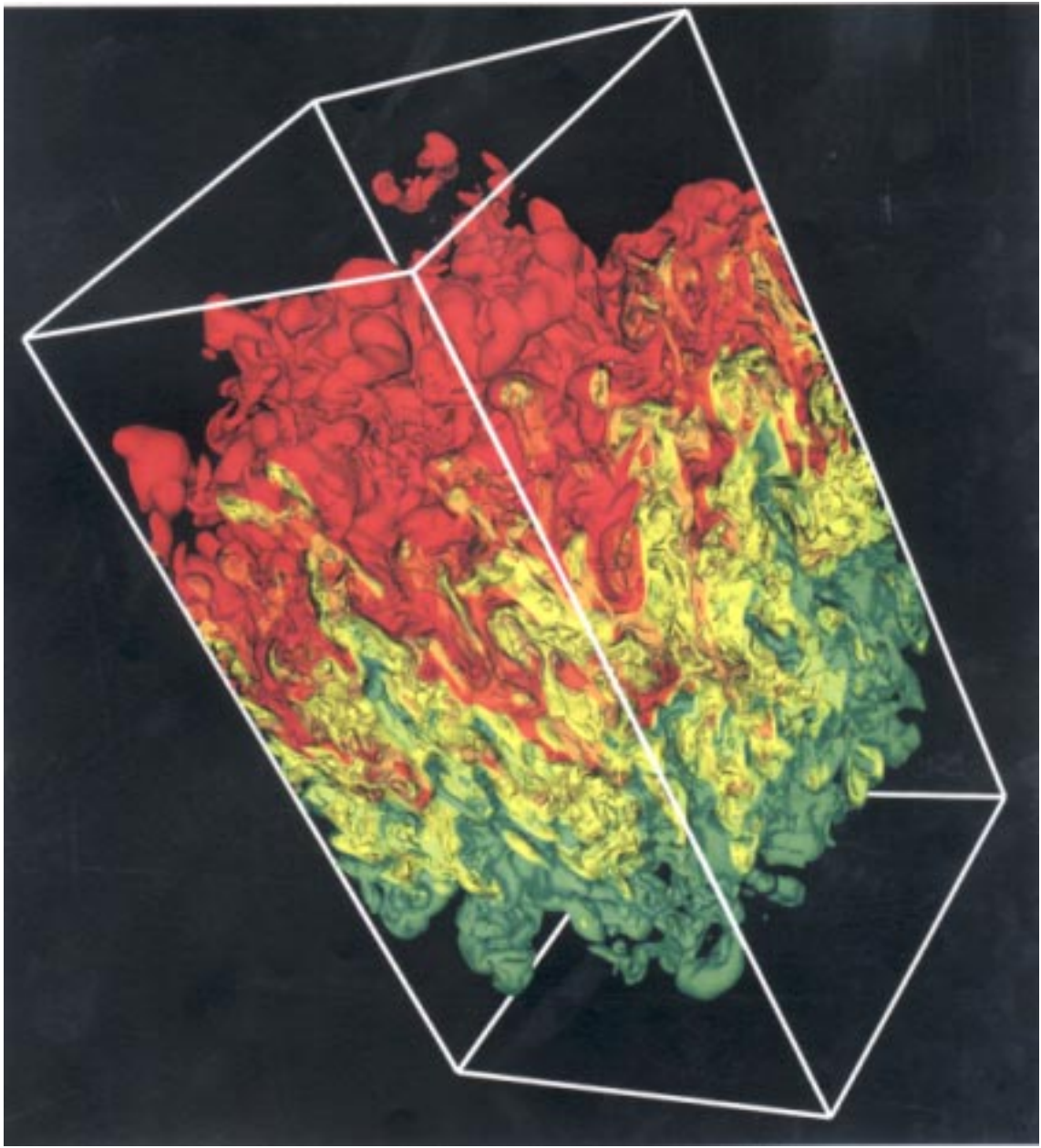

# Figure 7

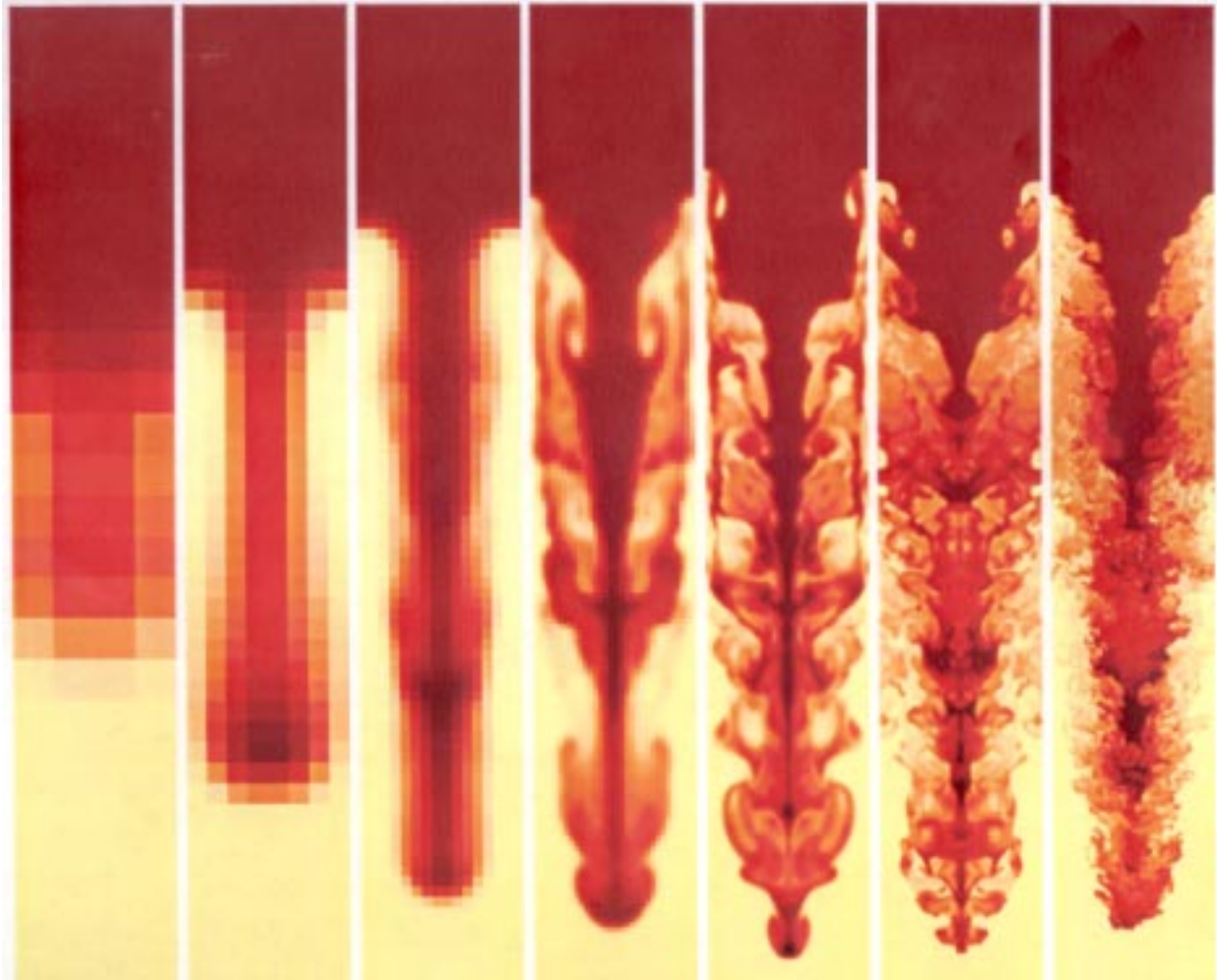

# Figure 8a

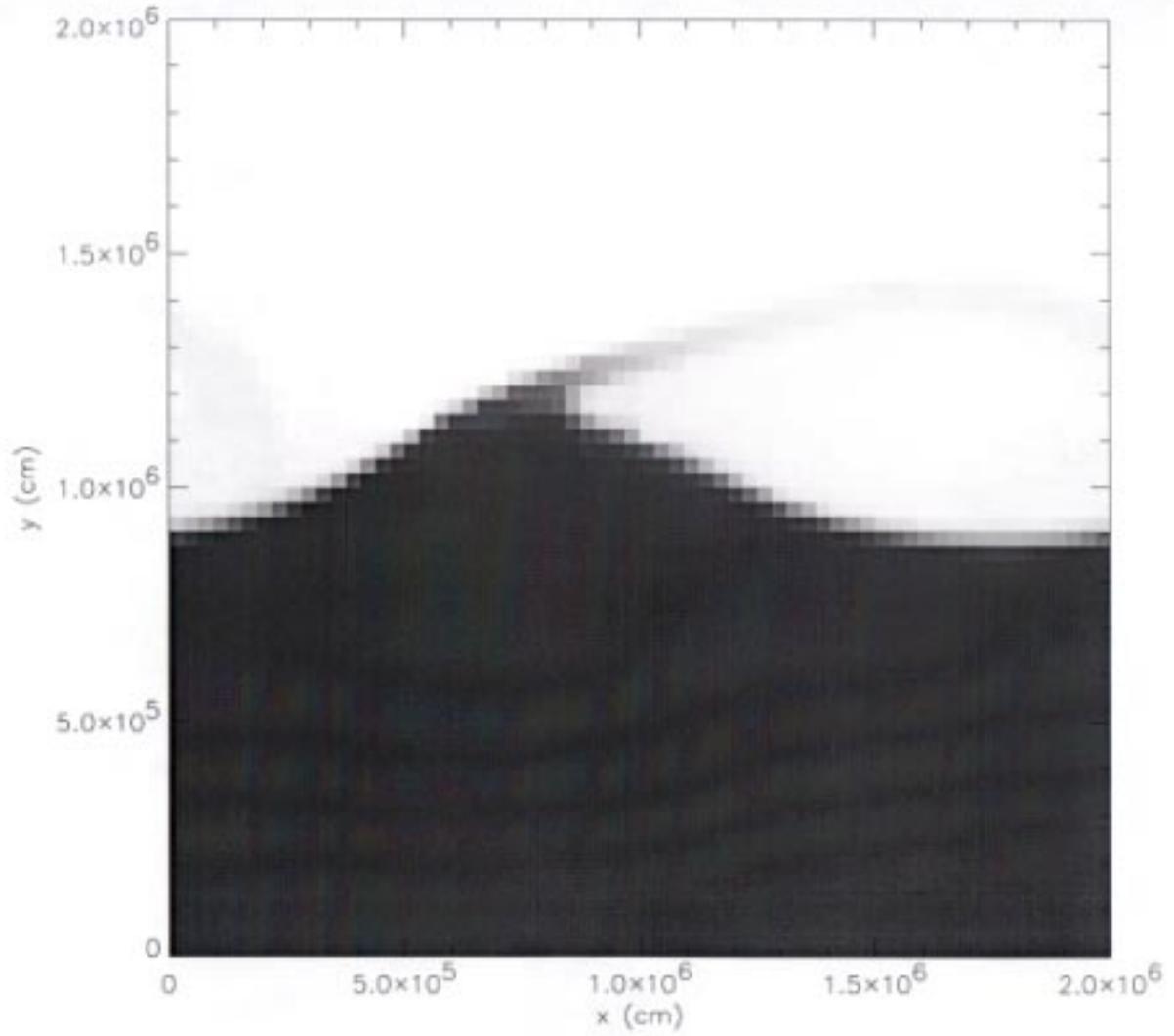

# Figure 8b

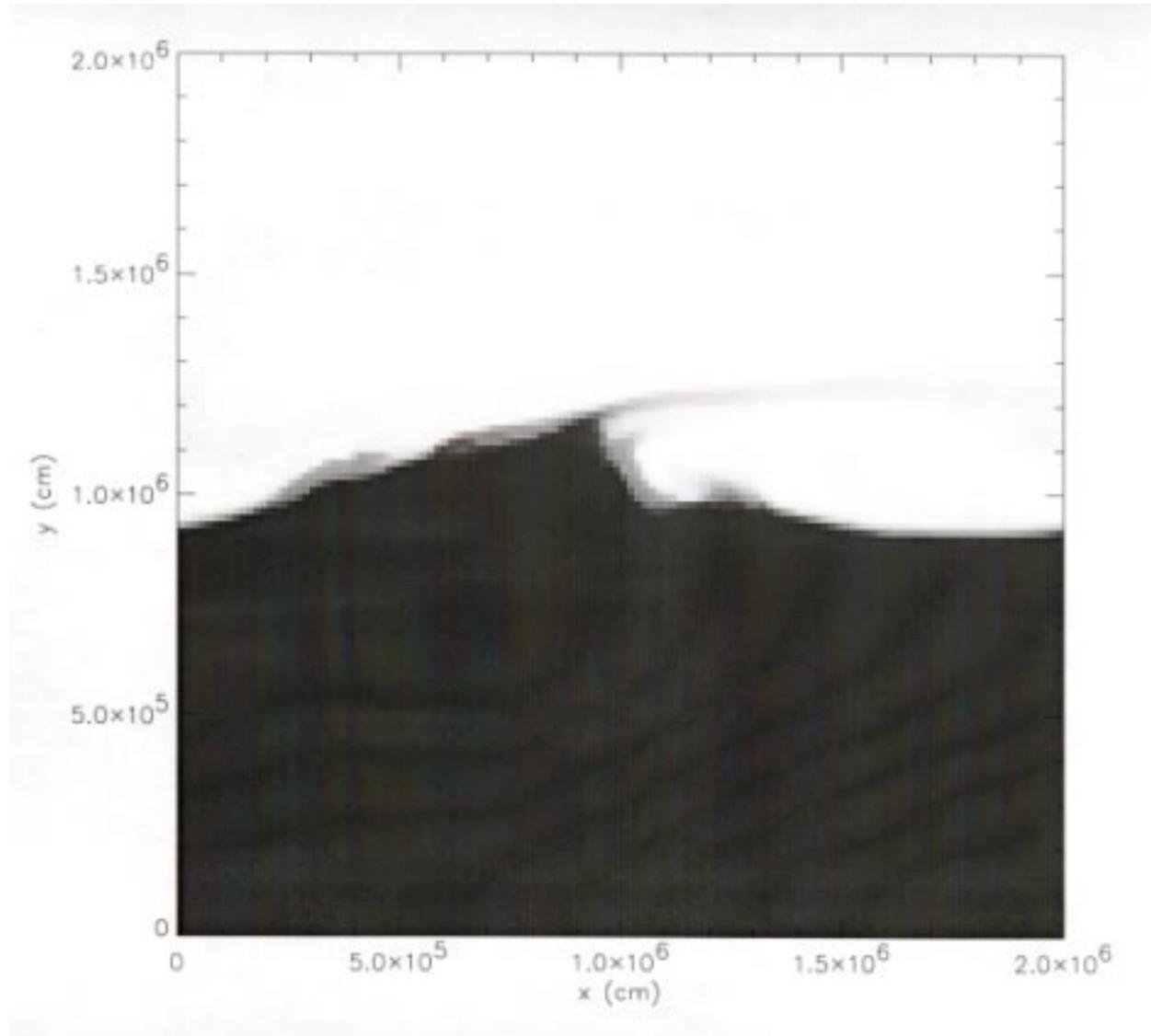

# Figure 8c

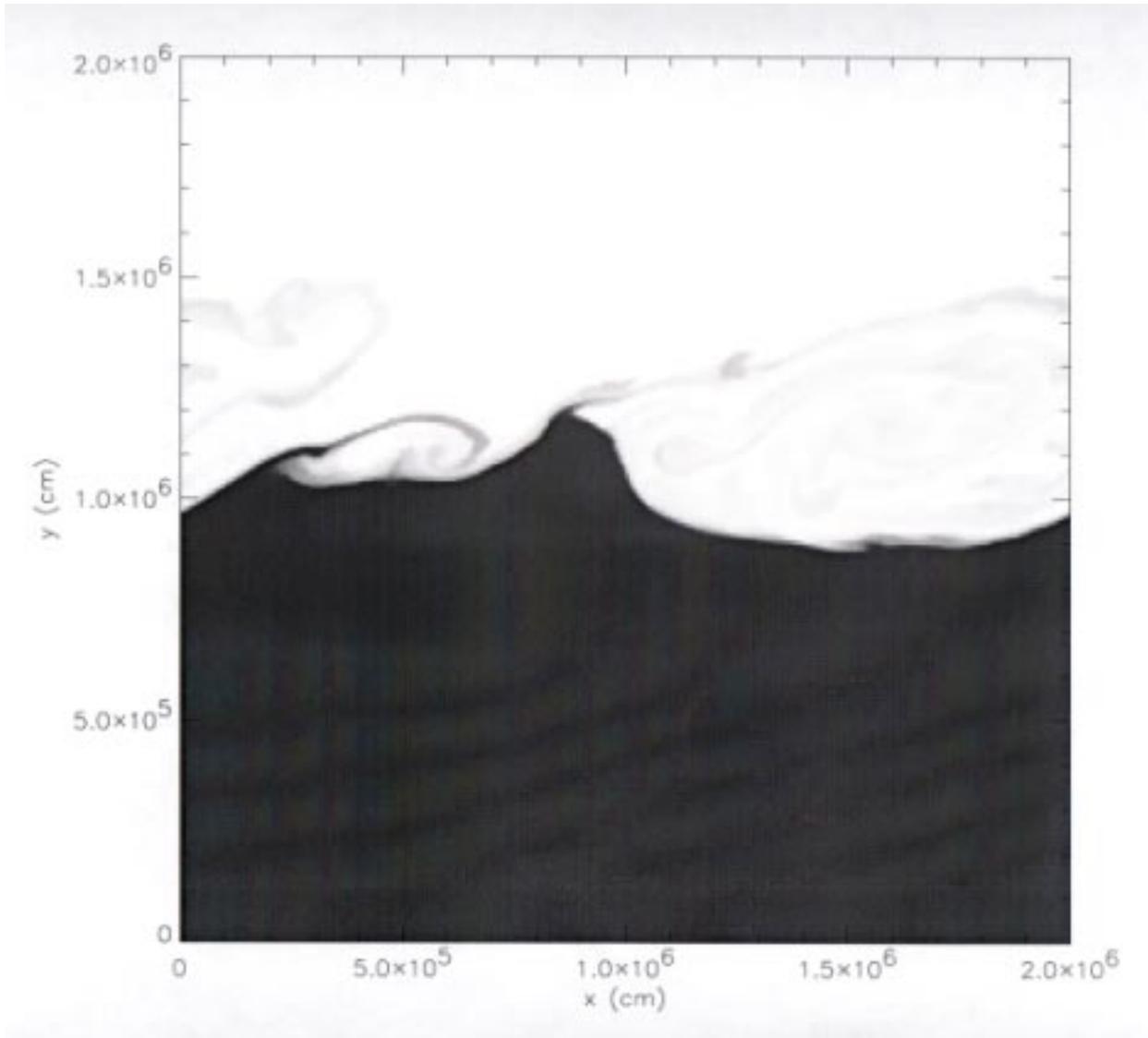

# Figure 8d

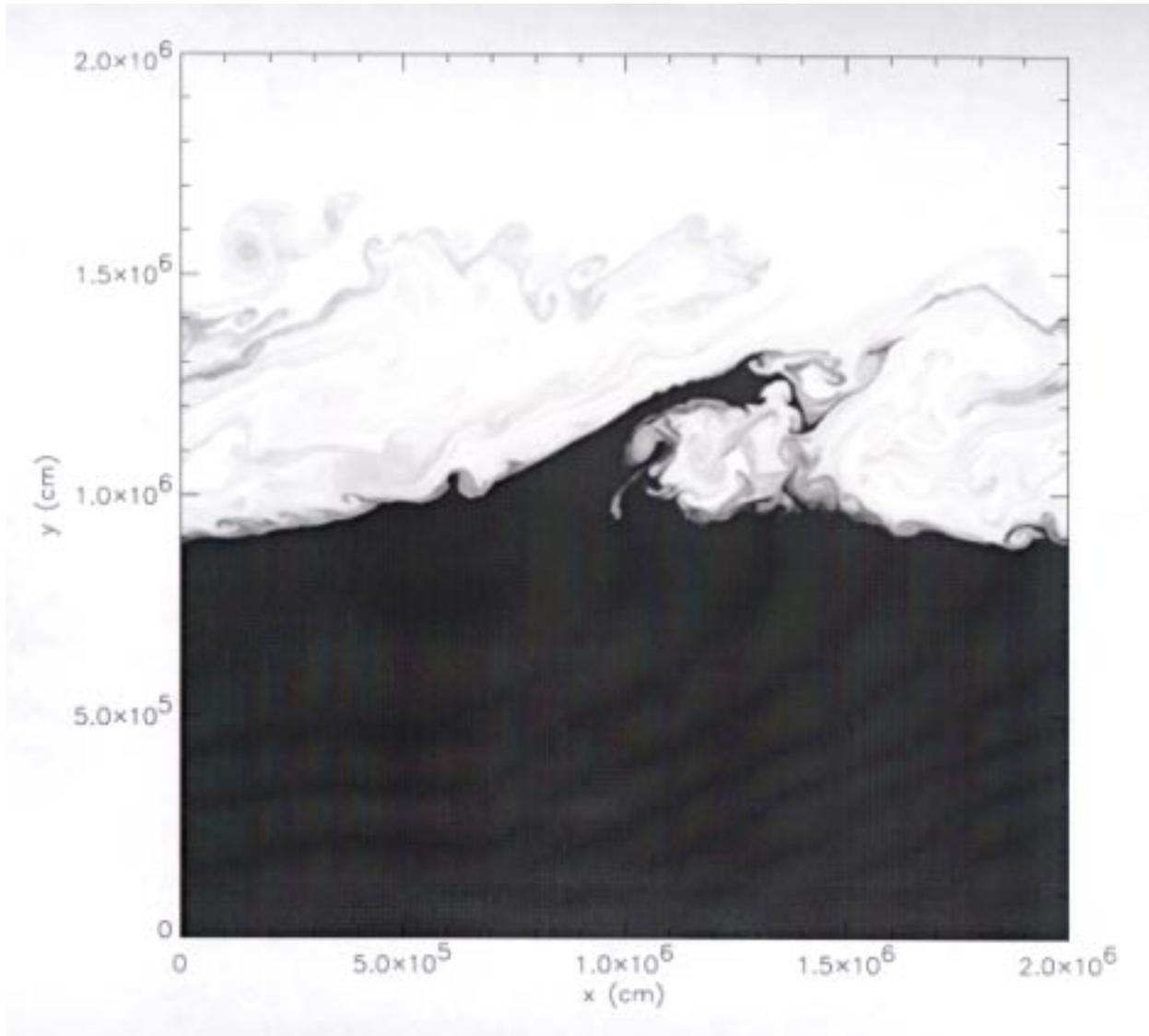

# Figure 8e

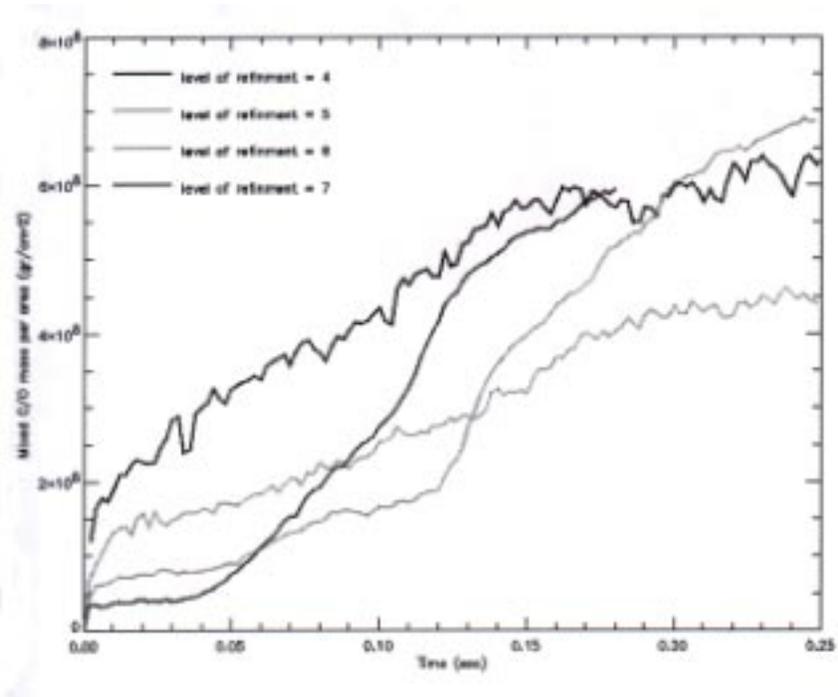